\documentclass{amsart}

\theoremstyle{definition}

\theoremstyle{remark}

\reversemarginpar

\begin{document}

\title[Elementary integrals]{The 
integrals in Gradshteyn and Ryzhik. \\
Part 7: Elementary examples}

\subjclass[2000]{Primary 33}


\author[Tewodros Amdeberhan]{Tewodros Amdeberhan}
\address{Department of Mathematics,
Tulane University, New Orleans, LA 70118}
\email{tamdeberhan@math.tulane.edu}

\author[Victor H. Moll]{Victor H. Moll}
\address{Department of Mathematics,
Tulane University, New Orleans, LA 70118}
\email{vhm@math.tulane.edu}

\thanks{The second  author wishes to the 
partial support of NSF-DMS 0409968.}

\begin{abstract}
The table of Gradshteyn and Ryzhik contains some elementary integrals. Some 
of them are derived here. 
\end{abstract}

\maketitle

\textwidth=6in
\textheight=8.5in

\newcommand{\nn}{\nonumber}
\newcommand{\ba}{\begin{eqnarray}}
\newcommand{\ea}{\end{eqnarray}}
\newcommand{\ift}{\int_{0}^{\infty}}
\newcommand{\no}{\noindent}
\newcommand{\at}{\text{tan}^{-1}}
\newcommand{\X}{{\mathbb{X}}}
\newcommand{\realpart}{\mathop{\rm Re}\nolimits}
\newcommand{\imagpart}{\mathop{\rm Im}\nolimits}

\newtheorem{Definition}{\bf Definition}[section]
\newtheorem{Thm}[Definition]{\bf Theorem} 
\newtheorem{Lem}[Definition]{\bf Lemma} 
\newtheorem{Cor}[Definition]{\bf Corollary} 
\newtheorem{Prop}[Definition]{\bf Proposition} 
\newtheorem{Note}[Definition]{\bf Note}
\numberwithin{equation}{section}

\maketitle

\section{Introduction} \label{intro} 
\setcounter{equation}{0}

Elementary mathematics leave the impression that 
there is marked difference between the two branches of calculus. 
{\em Differentiation} is a subject that is systematic: every evaluation
is a consequence of a number of rules and some basic examples.
However, {\em integration} is a mixture of art and science. The
successful evaluation 
of an integral depends on the right approach, the right change of variables or
a patient search in a table of integrals.  
In fact, the 
theory of {\em indefinite} integrals of elementary functions
is complete \cite{bronstein1,bronstein2}.  Risch's algorithm
determines whether a given function has an antiderivative 
within a given class of functions. 
However, the theory of {\em definite} integrals is far from complete and there
is no general theory available. The level of complexity in the evaluation of 
a definite integral is hard to predict as can be seen in
\begin{equation}
\ift e^{-x} \, dx = 1, \quad \ift e^{-x^{2}} \, dx = \frac{\sqrt{\pi}}{2}, 
\text{ and } \ift e^{-x^{3}} \, dx = \Gamma \left( \tfrac{4}{3} \right).  
\label{ex1}
\end{equation}
\noindent
The first integrand has an elementary primitive, the second one is 
the classical Gaussian integral and the evaluation of the third requires 
Euler's {\em gamma function} defined by
\begin{equation}
\Gamma(a) = \int_{0}^{\infty} x^{a-1}e^{-x} \, dx. 
\end{equation}

The table of integrals \cite{gr} contains a large variety of integrals. 
This 
paper continues the work initiated in \cite{moll-gr5,moll-gr1,
moll-gr2,moll-gr3, moll-gr4, moll-gr6} with 
the objective of providing proofs and 
context of {\em all the formulas} in the table \cite{gr}. Some
of them are truly elementary. In this paper we present a derivation of a 
small number of them. 

\section{A simple example} \label{sec-simple} 
\setcounter{equation}{0}

The first evaluation considered here is that of 
$\mathbf{3.249.6}$:
\begin{equation}
\int_{0}^{1} ( 1 - \sqrt{x})^{p-1} \, dx = \frac{2}{p(p+1)}.
\end{equation}
\noindent
The evaluation is completely elementary. The 
change of variables $y = 1 - \sqrt{x}$ produces 
\begin{equation}
I = -2 \int_{0}^{1} y^{p}\, dy + 2 \int_{0}^{1} y^{p-1} \, dy,
\end{equation}
\noindent
and each of these integrals can be evaluated directly to produce the 
result. 

This example can be generalized to consider 
\begin{equation}
I(a) = \int_{0}^{1} ( 1 - x^{a})^{p-1} \, dx.
\end{equation}
\noindent
The change of variables $t = x^{a}$ produces
\begin{equation}
I(a) = a^{-1} \int_{0}^{1} t^{1/a-1}(1-t)^{p-1} \, dt.
\end{equation}
\noindent
This integral appears as $\mathbf{3.251.1}$ and it 
can be  evaluated in terms of the 
{\em beta function} 
\begin{equation}
B(a,b) = \int_{0}^{1} x^{a-1}(1-x)^{b-1} \, dx,
\end{equation}
\noindent
as 
\begin{equation}
I(a)  = a^{-1} B \left( p, a^{-1} \right). 
\end{equation}
\noindent
The reader will find in \cite{moll-gr6} details about this evaluation.  \\

A further generalization is provided in the next lemma. \\

\begin{Lem}
Let $n \in \mathbb{N}, \, a, \, b, \, c \in \mathbb{R}$ with $bc > 0$. Define
$u = ac-b^{2}$. Then
\begin{equation}
\int_{0}^{1} \frac{a + b \sqrt{x}}{b + c \sqrt{x}} x^{n/2} \, dx = 
\frac{2u (-b)^{n+1}}{c^{n+3}} \ln(1 + c/b) + 
\frac{2u}{c^{2}} \sum_{j=0}^{n} \frac{(-1)^{j} }{n-j+1} 
\left( \frac{b}{c} \right)^{j} + \frac{2b}{(n+2)c}.
\nonumber
\end{equation}
\end{Lem}
\begin{proof}
Substitute $y = b + c \sqrt{x}$ and expand the new term $(y-b)^{n}$. 
\end{proof}

\section{A generalization of an algebraic example} \label{sec-gener} 
\setcounter{equation}{0}

The evaluation
\begin{equation}
\int_{-\infty}^{\infty} \frac{dx}{(1+x^{2}) \, \sqrt{4+3x^{2}} } = 
\frac{\pi}{3}
\end{equation}
\noindent
appears as $\mathbf{3.248.4}$ in \cite{gr}. We consider here the generalization
\begin{equation}
q(a,b) = \int_{-\infty}^{\infty} \frac{dx}{(1+x^{2}) \, \sqrt{b+ax^{2}} }.
\end{equation}
\noindent
We assume that $a, \, b > 0$. \\

The change of variables $x = \sqrt{b} t/\sqrt{a}$ yields 
\begin{equation}
q(a,b) = 2 \sqrt{a} \ift \frac{dt}{(a + bt^{2}) \, \sqrt{1+t^{2}} }
\end{equation}
\noindent
where we have used the symmetry of the integrand to write it over 
$[0, \infty)$. The standard trigonometric change of variables 
$t = \tan \varphi$ produces 
\begin{equation}
q(a,b) = 2 \sqrt{a} \int_{0}^{\pi/2} \frac{ \cos \varphi \, d \varphi}
{a \cos^{2}\varphi + b \sin^{2}\varphi}.
\end{equation}
\noindent
Finally, $u = \sin \varphi$, produces 
\begin{equation}
q(a,b) = 2 \sqrt{a} \int_{0}^{1} \frac{du}{a + (b-a)u^{2}}.
\end{equation}
\noindent
The evaluation of this integral is divided into three cases: \\

\noindent
{\bf Case 1}. $a=b$. Then we simply get $q(a,b) = 2/\sqrt{a}$. \\

\noindent
{\bf Case 2}. $a<b$. The change of variables $s = u \sqrt{b-a}/\sqrt{a}$
produces $(b-a)u^{2} = s^{2}a$, so that
\begin{equation}
q(a,b) = \frac{2}{\sqrt{b-a}} \int_{0}^{c} \frac{ds}{1+s^{2}} = 
\frac{2}{\sqrt{b-a}} \tan^{-1}c, 
\end{equation}
\noindent
with $c = \sqrt{b-a}/\sqrt{a}$. \\

\noindent
{\bf Case 3}. $a>b$. Then we write 
\begin{equation}
q(a,b) = 2 \sqrt{a} \int_{0}^{1} \frac{du}{a - (a-b)u^{2}}.
\end{equation}
\noindent
The change of variables $u = \sqrt{a} \, s/\sqrt{a-b}$ yields
\begin{equation}
q(a,b) = \frac{2}{\sqrt{a-b}} \int_{0}^{c} \frac{ds}{1-s^{2}},
\end{equation}
\noindent
where $c = \sqrt{a-b}/\sqrt{a}$. The partial fraction decomposition
\begin{equation}
\frac{1}{1-s^{2}} = \frac{1}{2} \left( \frac{1}{1+s} + \frac{1}{1-s} \right)
\end{equation}
\noindent
now produces 
\begin{equation}
q(a,b) = \frac{1}{\sqrt{a-b}} \, \ln \left( \frac{\sqrt{a} + \sqrt{a-b}}
{\sqrt{a}-\sqrt{a-b}} \right).
\end{equation}
\noindent
The special case in $\mathbf{3.248.4}$ 
corresponds to $a=3$ and $b=4$. The value of the 
integral is $2 \tan^{-1}(1/\sqrt{3}) = \frac{\pi}{3}$, as claimed.  This 
generalization has been included as ${\mathbf{3.248.6}}$ in 
the latest edition of \cite{gr}. \\

We now consider a generalization of this integral. The proof requires several
elementary steps, given first for the convenience of the  reader.  \\

Let $a, \, b \in \mathbb{R}$ with $a < b$ and $n \in \mathbb{N}$.
Introduce the notation 
\begin{equation}
I = I_{n}(a,b) := 
\int_{0}^{\infty} \frac{dt}{(a+bt^{2})^{n} \, \sqrt{1+t^{2}}}. 
\end{equation}
\noindent
Then we have: \\

\begin{Lem}
The integral $I_{n}(a,b)$ is given by 
\begin{equation}
I_{n}(a,b) = \int_{0}^{1} \frac{(1-v^{2})^{n-1} \, dv}{(a + \alpha v^{2})^{n}}, 
\end{equation}
\noindent
with $\alpha = b -a$. 
\end{Lem}
\begin{proof}
The change of variables $v = t/\sqrt{1+t^{2}}$ gives the result. 
\end{proof}

The identity 
\begin{equation}
(1-v^{2})^{n-1} = (1-v^{2})^{n} + (1-v^{2})^{n-1} \left\{ \frac{1}{\alpha} 
(a + \alpha v^{2}) - \frac{a}{\alpha} \right\}
\end{equation}
\noindent
produces 
\begin{equation}
I_{n}(a,b) = \frac{\alpha}{b} \int_{0}^{1} \frac{(1-v^{2})^{n}}{(a + 
\alpha v^{2})^{n}} \, dv + \frac{1}{b} \int_{0}^{1} 
\frac{(1-v^{2})^{n-1} }{(a + \alpha v^{2})^{n-1}} \, dv. 
\end{equation}

The evaluation of these integrals requieres an intermediate result, that is 
also of independent interest. 

\begin{Lem}
\label{lemma31}
Assume $z \in \mathbb{R}$ and $n \in \mathbb{N} \cup \{ 0 \}$. Then 
\begin{equation}
\int_{0}^{1} \frac{dx}{(1+ z^{2}x^{2})^{n+1} } = 
\frac{1}{2^{2n}} \binom{2n}{n} \left( \frac{\tan^{-1}z}{z} + 
\sum_{k=1}^{n} \frac{2^{2k}}{2k \binom{2k}{k}} \frac{1}{(1+z^{2})^{k}} 
\right). 
\end{equation}
\end{Lem}
\begin{proof}
Define 
\begin{equation}
F_{n}(z) := \int_{0}^{1} \frac{dx}{(1+z^{2}x^{2})^{n+1}} = 
\frac{1}{z} \int_{0}^{z} \frac{dy}{(1+y^{2})^{n+1}}. 
\label{iden1}
\end{equation}
\noindent
Take derivatives with respect to $z$ on both sides of (\ref{iden1}).  The 
outcome is a system of differential-difference equations
\begin{eqnarray}
\frac{dF_{n}(z)}{dz}  & = & \frac{2(n+1)}{z} F_{n+1}(z) - \frac{2(n+1)}{z} 
F_{n}(z) \label{system1} \\
\frac{dF_{n}}{dz} & = & - \frac{1}{z} F_{n}(z) + \frac{1}{z(1+z^{2})^{n+1}}. 
\nonumber 
\end{eqnarray}
\noindent
Solving for a purely recursive relation we obtain (after re-indexing 
$n \mapsto n-1$):
\begin{equation}
2nF_{n}(z) = (2n-1)F_{n-1}(z) + \frac{1}{(1+z^{2})^{n}},
\end{equation}
\noindent
with the initial condition $F_{0}(z) = \frac{1}{z} \tan^{-1}z$. This recursion
is solved using the procedure described in Lemma $2.7$ of 
\cite{moll-gr5}. This produces the stated expression for $F_{n}(z)$. 
\end{proof}

The next required evaluation is that of the powers of a simple rational 
function.

\begin{Lem}
Let $a, \, b, \, c, \, d$ be real numbers such that $cd > 0$. Then 
\begin{eqnarray}
\int_{0}^{1} \left( \frac{ax^{2} + b}{cx^{2}+d} \right)^{n} \, dx  & = & 
\frac{a^{n}}{c^{n}} + \frac{4a^{n}}{c^{n}} \sqrt{ \frac{d}{c}} 
\tan^{-1} \sqrt{c/d} \sum_{k=1}^{n} \left( \frac{bc-ad}{4ad} \right)^{k} 
\binom{n}{k} \binom{2k-2}{k-1} \nonumber \\
& + & \frac{4a^{n}}{c^{n}} \sum_{k=1}^{n} \left( \frac{bc-ad}{4ad} \right)^{k} 
\binom{n}{k} \binom{2k-2}{k-1} 
\sum_{j=1}^{k-1} \frac{2^{2j}}{2j \, \binom{2j}{j} } 
\frac{d^{j}}{(c+d)^{j}}. 
\nonumber 
\end{eqnarray} 
\end{Lem}
\begin{proof}
Start with the partial fraction expansion 
\begin{equation}
G(x) := \frac{ax^{2}+b}{cx^{2}+d} = \frac{a}{c} + 
\frac{bc-ad}{cd} \frac{1}{cx^{2}/d + 1},
\end{equation}
\noindent
and expand  $G(x)^{n}$ by the binomial theorem. The result follows by
using Lemma \ref{lemma31}. 
\end{proof}

The next result follows by combining the statements of the previous three 
lemmas.

\begin{Thm}
Let $a, \, b \in {\mathbb{R}}^{+}$ with $a< b$. Then
\begin{eqnarray}
I_{n+1}(a,b) & := & \int_{0}^{\infty} \frac{dt}{(a+bt^{2})^{n+1} \, \sqrt{1+t^{2}} }   \nonumber \\
& = &  \frac{1}{a(a-b)^{n}} \sum_{j=0}^{n} \binom{n}{j} 
\left( \frac{-b}{4a} \right)^{j} \binom{2j}{j} \times
\left( 
\frac{\tan^{-1} \sqrt{b/a-1} }{\sqrt{b/a-1}} + 
\sum_{k=1}^{j} \frac{2^{2k}}{2k \binom{2k}{k} } \left( \frac{a}{b} \right)^{k} 
\right). 
\nonumber
\end{eqnarray}
\end{Thm}

\section{Some integrals involving the exponential function} \label{sec-simexp} 
\setcounter{equation}{0}

In \cite{gr} we find $\mathbf{3.310}$: 
\begin{equation}
\ift e^{-px} \, dx = \frac{1}{p}, \text{   for } p >0,
\end{equation}
\noindent
that is probably the most elementary evaluation in the table. The example 
$\mathbf{3.311.1}$
\begin{equation}
\ift \frac{dx}{1+ e^{px}} = \frac{\ln 2}{p},
\label{int-99}
\end{equation}

\noindent
can also be evaluated in elementary terms. Observe first that the change of 
variables $t = px$, shows that (\ref{int-99}) is equivalent to the case $p=1$:
\begin{equation}
\ift \frac{dx}{1+ e^{x}} = \ln 2.
\label{int-99a}
\end{equation}
\noindent
This can be evaluated by the change of variables $u = e^{x}$ that yields
\begin{equation}
I = \int_{1}^{\infty} \frac{du}{u(1+u)},
\end{equation}
\noindent
and this can be integrated by partial fractions to produce the result. The 
parameter in (\ref{int-99}) is what we call {\em fake}, in the sense that
the corresponding integral is independent of it. The 
advantage of such  parameter is that it provides flexibility to a formula: 
differentiating (\ref{int-99}) with respect to $p$ produces 
\begin{equation}
\ift \frac{xe^{px} \, dx}{(1+ e^{px})^{2}} = \frac{\ln 2}{p^{2}}, 
\end{equation}
\begin{equation}
\ift \frac{x^{2}e^{px}(e^{px}-1) \, dx}{(1+ e^{px})^{3}} = 
\frac{2 \ln 2}{p^{3}}, 
\end{equation}
\begin{equation}
\ift \frac{x^{3}e^{px}(e^{2px}-4e^{px}+1) \, dx}{(1+ e^{px})^{4}} = 
\frac{6 \ln 2}{p^{4}}. 
\end{equation}
\noindent
The general integral formula  is obtained by differentiating (\ref{int-99}) 
$n$-times with respect to $p$ to produce
\begin{equation}
\ift \left( \frac{\partial}{\partial p} \right)^{n} \frac{dx}{1+e^{px}} = 
(-1)^{n} \frac{n!}{p^{n+1}} \ln 2. 
\end{equation}
\noindent
The pattern of the integrand is clear:
\begin{equation}
\left( \frac{\partial }{\partial p} \right)^{n} \frac{1}{1+e^{px}} = 
\frac{(-1)^{n}x^{n} e^{px} }{(1+e^{px})^{n+1}} P_{n}(e^{px}),
\end{equation}
\noindent
where $P_{n}$ is a polynomial of degree $n-1$. It follows that 
\begin{equation}
\int_{0}^{\infty} \frac{x^{n} e^{px} P_{n}(e^{px}) \, dx }
{(1+e^{px})^{n+1}} = \frac{n! \, \ln 2}{p^{n+1}}.
\end{equation}
\noindent
The change of variables $t = px$ shows that $p$ is a fake parameter. The 
integral is equivalent to
\begin{equation}
\int_{0}^{\infty} \frac{x^{n} e^{x} P_{n}(e^{x}) \, dx }
{(1+e^{x})^{n+1}} = n! \, \ln 2.
\end{equation}

The first few polynomials in the sequence are given by
\begin{eqnarray}
P_{1}(u) & = & 1, \\
P_{2}(u) & = & u-1, \nonumber \\
P_{3}(u) & = & u^{2}-4u+1, \nonumber \\
P_{4}(u) & = & u^{3}-11u^{2}+11u-1. \nonumber
\end{eqnarray}

\medskip

\begin{Prop}
The polynomials $P_{n}(u)$ satisfy the recurrence 
\begin{equation}
P_{n+1}(u) = (nu-1)P_{n}(u) -u(1+u) \frac{d}{du} P_{n}(u).
\label{recP}
\end{equation}
\end{Prop}
\begin{proof}
The result follows by expanding the relation
\begin{equation}
\frac{(-1)^{n+1}x^{n+1} e^{px} P_{n+1}(e^{px}) }{(1+e^{px})^{n+2}} = 
\frac{\partial}{\partial p} \left( \frac{(-1)^{n}x^{n} e^{px} P_{n}(e^{px}) }
{(1+e^{px})^{n+1}} \right).
\end{equation}
\end{proof}

\medskip

Examining the first few values, we observe that 
\begin{equation}
Q_{n}(u) := (-1)^{n} P_{n}(-u)
\end{equation}
\noindent
is a polynomial with positive coefficients. This follows directly from the 
recurrence
\begin{equation}
Q_{n+1}(u) = (nu+1)Q_{n}(u) +u(1-u) \frac{d}{du} Q_{n}(u).
\end{equation}
\noindent
This comes directly from (\ref{recP}).  The first few values are
\begin{eqnarray}
Q_{1}(u) & = & 1, \\
Q_{2}(u) & = & u+1, \nonumber \\
Q_{3}(u) & = & u^{2}+4u+1, \nonumber \\
Q_{4}(u) & = & u^{3} + 11u^{2}+11u+1. \nonumber
\end{eqnarray}

Writing 
\begin{equation}
Q_{n}(u) = \sum_{j=0}^{n-1} E_{j,n}u^{j},
\end{equation}
\noindent
the reader will easily verify the recurrence
\begin{eqnarray}
E_{0,n+1} & = & E_{0,n} \label{rec-Euler1} \\ 
E_{j,n+1} & = & (n-j+1)E_{j-1,n} + (j+1)E_{j,n} \nonumber \\
E_{n,n+1} & = & E_{n,n}. \nonumber
\end{eqnarray}
\noindent
The numbers $E_{j,n}$ are called {\em Eulerian numbers}. They appear in 
many situations. For example, they provide the coefficients in the series
\begin{equation}
\sum_{k=1}^{\infty} k^{j}x^{k} = \frac{x}{(1-x)^{j+1}} 
\sum_{n=0}^{m-1} E_{j,n}x^{n}
\end{equation}
\noindent
and satisfy the simpler recurrence 
\begin{equation}
E_{j,n} = nE_{j-1,n} + jE_{j,n-1},
\end{equation}
\noindent
that can be derived from (\ref{rec-Euler1}). These numbers have a combinatorial
interpretation: they count the number of permutations of $\{ 1, \, 2, \ldots, 
n \}$ having $j$ permutation ascents. The explicit formula
\begin{equation}
E_{j,n} = \sum_{k=0}^{j+1} (-1)^{k} \binom{n+1}{k} (j+1-k)^{n},
\end{equation}
\noindent
can be checked from the recurrences.  The reader will find more information
about these numbers in \cite{graham1}.

\section{A simple change of variables} \label{sec-simcha} 
\setcounter{equation}{0}

The table \cite{gr} contains the example $\mathbf{3.195}$:
\begin{equation}
\ift \frac{(1+x)^{p-1}}{(x+a)^{p+1}} \, dx = \frac{1-a^{-p}}{p(a-1)}.
\label{int-simple}
\end{equation}
\noindent
One must include the restrictions $a>0, \, a \neq 1, \, p \neq 0$. The 
evaluation is elementary: let 
\begin{equation}
u = \frac{1+x}{x+a},
\end{equation}
\noindent
to obtain 
\begin{equation}
I = \frac{1}{a-1} \int_{1/a}^{1} u^{p-1} \, du,
\label{int-100}
\end{equation}
\noindent
that gives the stated value. The formula has been supplemented with the 
value $1$ for $a=1$ and $\ln a/(a-1)$ when $p=0$ in the last edition 
of \cite{gr}. \\

Differentiating (\ref{int-simple}) $n$ times with respect to the parameter $p$
produces
\begin{equation}
\ift \frac{(1+x)^{p-1}}{(x+a)^{p+1}} \ln^{n} \left( \frac{1+x}{x+a} \right) \,
dx = \frac{(-1)^{n} a^{-p} }{(a-1)p^{n+1}} 
\left[ n! \, (a^{p}-1) - \sum_{k=1}^{n} \frac{n! (p \, \ln a)^{k}}{k!} \right].
\nonumber
\end{equation}
\noindent
Naturally, the integral above is just 
\begin{equation}
\frac{1}{a-1} \int_{1/a}^{1} u^{p-1} \, \ln^{n}u \, du
\end{equation}
\noindent
and its value can also be obtained by differentiation of (\ref{int-100}).  \\

The next result presents a generalization of (\ref{int-simple}):

\begin{Lem}
\label{free1}
Let $a, \, b$ be free parameters and $n \in \mathbb{N}$. Then
\begin{equation}
\int_{0}^{\infty} \frac{(1+x)^{b}}{(x+a)^{b+n}} \, dx
= (a-1)^{-n} \times \left\{ B(n,b) -  
\sum_{k=0}^{n-1} (-1)^{k} \binom{n-1}{k} \frac{a^{-b-k}}{b+k}
\right\}, 
\nonumber
\end{equation}
\noindent
where $B(n,b)$ is Euler's beta function.
\end{Lem}
\begin{proof}
Use the change of variables $u=(1+x)/(a+x)$, expand in series and then integrate
term by term.
\end{proof}

\section{Another example} \label{sec-another}
\setcounter{equation}{0}

The integral in $\mathbf{3.268.1}$ states that 
\begin{equation}
\int_{0}^{1} \left( \frac{1}{1-x} - \frac{px^{p-1}}{1-x^{p}} \right) \, dx 
= \ln p.
\end{equation}
\noindent
To compute it, and to avoid the singularity at $x=1$, we write
\begin{equation}
I = \lim\limits_{\epsilon \to 0} \int_{0}^{1 - \epsilon} 
 \left( \frac{1}{1-x} - \frac{px^{p-1}}{1-x^{p}} \right) \, dx.
\end{equation}
\noindent
This evaluates as 
\begin{equation}
I = \lim\limits_{\epsilon \to 0} -\ln \epsilon + 
\ln( 1 - (1- \epsilon)^{p}) = 
\lim\limits_{\epsilon \to 0} \, \ln \left( \frac{1 - (1 - \epsilon)^{p}}
{\epsilon} \right) = \ln p. 
\end{equation}

\section{Examples of recurrences} \label{sec-recurr}
\setcounter{equation}{0}

Several definite integrals in \cite{gr} can be evaluated by producing a 
recurrence for them. For example, in $\mathbf{3.622.3}$ we find 
\begin{equation}
\int_{0}^{\pi/4} \tan^{2n}x \, dx = (-1)^{n} 
\left( \frac{\pi}{4} - \sum_{j=0}^{n-1} \frac{(-1)^{j-1}}{2j-1} \right). 
\label{tan-rec}
\end{equation}
\noindent
To check this identity, define 
\begin{equation}
I_{n} = \int_{0}^{\pi/4} \tan^{2n}x \, dx 
\end{equation}
\noindent 
and integrate by parts to produce 
\begin{equation}
I_{n} = -I_{n-1} + \frac{1}{2n-1}.
\label{tan-recur}
\end{equation}
\noindent
From here we generate the first few values 
\begin{equation}
I_{0} = \frac{\pi}{4}, \, 
I_{1} = -\frac{\pi}{4}+ 1, \, 
I_{2} = \frac{\pi}{4}- 1+ \frac{1}{3}, \text{ and }
I_{3} = -\frac{\pi}{4}+ 1- \frac{1}{3}+ \frac{1}{5},  \nonumber
\end{equation}
\noindent
and from here one can {\em guess} the formula (\ref{tan-rec}). A proof by 
induction is easy using (\ref{tan-recur}). \\

A similar argument produces $\mathbf{3.622.4}$:
\begin{equation}
\int_{0}^{\pi/4} \tan^{2n+1}x \, dx =  \frac{(-1)^{n+1}}{2}
\left( \ln 2 - \sum_{k=1}^{n} \frac{(-1)^{k}}{k} \right). 
\label{tan-recJ}
\end{equation}
\noindent
To establish this, define
\begin{equation}
J_{n} = \int_{0}^{\pi/4} \tan^{2n+1}x \, dx
\end{equation}
\noindent
and integrate by parts to produce 
\begin{equation}
J_{n}= - J_{n-1} + \frac{1}{2n}.
\label{rec-J}
\end{equation}
\noindent
The value 
\begin{equation}
J_{0} = \int_{0}^{\pi/4} \tan x \, dx = \frac{\ln 2}{2},
\end{equation}
\noindent
and the recurrence (\ref{rec-J}) yield the formula.

\section{A truly elementary example} \label{sec-truly}
\setcounter{equation}{0}

The evaluation of $\mathbf{3.471.1}$
\begin{equation}
\int_{0}^{u} \text{exp} \left( - \frac{b}{x} \right) \frac{dx}{x^{2}} 
= \frac{1}{b} \text{exp} \left( - \frac{b}{u} \right),
\end{equation}
\noindent
is truly elementary: the change of variables $t = -b/x$ gives the result.

\section{Combination of polynomials and exponentials}
\setcounter{equation}{0}

Integration by parts produces
\begin{equation}
\int x^{n} e^{ax} \,dx = \frac{1}{a} x^{n}e^{ax} -
\frac{n}{a} \int x^{n-1} e^{ax}.
\label{rec-exp}
\end{equation}
\noindent
This appears as $\mathbf{2.321.1}$ in \cite{gr}. Introduce the notation
\begin{equation}
I_{n}(a) := \int x^{n}e^{ax} \, dx
\end{equation}
\noindent
so that (\ref{rec-exp}) states that
\begin{equation}
I_{n}(a) = \frac{1}{a} x^{n} e^{ax} - \frac{n}{a} I_{n-1}(a).
\end{equation}
\noindent
This recurrence is now used to prove
\begin{equation}
I_{n}(a) = n! e^{ax} \sum_{k=0}^{n} \frac{(-1)^{k} x^{n-k}}{(n-k)! \, 
a^{k+1}}
\label{int-rec11}
\end{equation}
\noindent
by an easy induction argument. This appears as $\mathbf{2.321.2}$. The 
case $ 1 \leq n \leq 4$ appear as $\mathbf{2.322.1}, \mathbf{2.322.2}, 
\mathbf{2.322.3}, \mathbf{2.322.4}$, respectively.  \\

Integrating (\ref{int-rec11}) between $0$ and $u$ produces $\mathbf{3.351.1}$:
\begin{equation}
\int_{0}^{u} x^{m} e^{-ax} \, dx = \frac{n!}{a^{n+1}} - 
e^{-au} \sum_{k=0}^{n} \frac{n!}{k!} \frac{u^{k}}{a^{n-k+1}}.
\end{equation}
\noindent
This sum can be written in terms of the incomplete gamma function. Details 
will appear in a future publication. Integrating
(\ref{int-rec11}) from $u$ to 
$\infty$ produces
\begin{equation}
\int_{u}^{\infty}x^{n}e^{-ax} dx = e^{-au} \sum_{k=0}^{n} 
\frac{n!}{k!} \frac{u^{k}}{a^{n-k+1}}. 
\end{equation}
\noindent
This appears as $\mathbf{3.351.2}$. \\

The special case $n=1$ of $3.351.1$ appears as 
$\mathbf{3.351.7}$. The cases $n=2$ and $n=3$ appear as $\mathbf{3.351.8}$
and $\mathbf{3.351.9}$, respectively.

\section{A perfect derivative} \label{sec-perfect} 
\setcounter{equation}{0}

In section $\mathbf{4.212}$ we find a list of examples that can be evaluated in 
terms of the exponential integral function, defined by
\begin{equation}
\text{Ei}(x) := \int_{-\infty}^{x} \frac{e^{t} \, dt}{t}
\label{eidef}
\end{equation}
\noindent
for $x<0$ and by the Cauchy principal value of (\ref{eidef}) for $x > 0$.
An exception is 
$\mathbf{4.212.7}$:
\begin{equation}
\int_{1}^{e} \frac{\ln x \, dx}{(1+ \ln x)^{2}} = \frac{e}{2} - 1. 
\label{int-simexp}
\end{equation}

\noindent 
This is an elementary integral: the change of variables $t = 1 + \ln x$ 
yields
\begin{equation}
I = \frac{1}{e} \int_{1}^{2} \frac{(t-1)}{t^{2}}e^{t} \, dt
\end{equation}
\noindent
and to evaluate it observe that 
\begin{equation}
\frac{(t-1)}{t^{2}}e^{t} = \frac{d}{dt} \frac{e^{t}}{t}. 
\end{equation}

The change of variables $t = \ln x$ in (\ref{int-simexp}) yields
\begin{equation}
\int_{0}^{1} \frac{t \, e^{t} \, dt}{(1+t)^{2}} = \frac{e}{2}-1.
\end{equation}

\noindent
This is $\mathbf{3.353.4}$ in \cite{gr}. \\

The previous evaluation can be generalized by introducing a parameter. 

\begin{Lem}
Let $\alpha \in \mathbb{R}$. Then
\begin{equation}
\int_{1}^{e} \frac{\ln x \, dx }{(\alpha + \ln x)^{\alpha+1}} = 
\frac{e}{(\alpha+1)^{\alpha}} - \frac{1}{\alpha^{\alpha}}.
\end{equation}
\end{Lem}
\begin{proof}
Substitute $t = \alpha + \ln x$ and use 
\begin{equation}
\frac{d}{dt} \frac{e^{t}}{t^{\alpha}} = \frac{t-\alpha}{t^{\alpha+1}} e^{t}.
\end{equation}
\noindent
The case $\alpha = 1$ corresponds to (\ref{int-simexp}). 
\end{proof}

\section{Integrals involving quadratic polynomials} \label{sec-quadratic} 
\setcounter{equation}{0}

There are several evaluation in \cite{gr} that involve quadratic polynomials. 
We assume, for reasons of convergence, that the discriminant 
$d = b^{2}-ac$ is strictly negative. \\

We start with 
\begin{equation}
\ift \frac{dx}{ax^{2}+2bx+c} = \frac{1}{\sqrt{ac-b^{2}}}
\cot^{-1} \left( \frac{b}{\sqrt{ac-b^{2}}} \right).
\label{base-case}
\end{equation}
\noindent
This is evaluated by completing the square and a simple trigonometric 
substitution:
\begin{eqnarray}
\ift \frac{dx}{ax^{2}+2bx+c} & = &  \frac{1}{a} \int_{b/a}^{\infty} 
\frac{du}{u^{2} -d/a^{2}} \nonumber \\
& = & \frac{1}{\sqrt{-d}} \int_{b/\sqrt{-d}}^{\infty} \frac{dv}{v^{2}+1}.
\nonumber 
\end{eqnarray}
\noindent
Differentiating (\ref{base-case})  with respect to $c$ produces 
$\mathbf{3.252.1}$:

\begin{equation}
\ift \frac{dx}{(ax^{2}+2bx+c)^{n}} = \frac{(-1)^{n-1}}{(n-1)!} 
\frac{\partial^{n-1}}{\partial c^{n-1}} 
\left[ \frac{\cot^{-1}(b/\sqrt{ac-b^{2}})}{\sqrt{ac-b^{2}}} \right]. 
\nonumber
\end{equation}

We now produce a closed-from expression for this integral.  \\

\begin{Lem}
\label{ex2}
Let $n \in \mathbb{N}$ and $u := 4(ac-b^{2})/ac$. Assume $cu>0$. Then we have
the explicit evaluation
\begin{equation}
\ift \frac{dx}{(ax^{2}+2bx+c)^{n}} = 
\frac{2b}{a(cu)^{n}} \binom{2n-2}{n-1} \times 
\left\{ \frac{\sqrt{acu}}{b} \cot^{-1} \left( \frac{2b}{\sqrt{acu}} \right) -
\sum_{j=1}^{n-1} \frac{u^{j}}{j \binom{2j}{j}}
\right\}.
\label{explicit1}
\end{equation}
\end{Lem}
\begin{proof}
The case $n=1$ was described above:
\begin{equation}
h(a,b,c) := \int_{0}^{\infty} \frac{dx}{ax^{2}+2bx + c} = 
\frac{1}{\sqrt{ac-b^{2}}} \cot^{-1} \left( \frac{1}{\sqrt{ac-b^{2}}} \right).
\end{equation}
\noindent
Now observe that $h(a^{2},abc,b^{2}) = h(1,b,1)/ac$. To establish 
(\ref{explicit1}), change the parameters sequentially as $a \mapsto a^{2}; \, 
c \mapsto c^{2}; \, b \mapsto abc$. In the new format, both sides satisfy
the differential-difference equation
\begin{equation}
-2nc(1-b^{2}) f_{n+1} = \frac{df_{n}}{dc} + \frac{b}{ac^{2n}}. 
\end{equation}
\noindent
The identity (\ref{explicit1}) is obtained by reversing the transformations
of paramaters indicated above. 
\end{proof}

\medskip

\begin{Cor}
Using the notations of Lemma \ref{ex2} we have
\begin{equation}
\sum_{j=1}^{\infty} \frac{u^{j}}{j \binom{2j}{j}} = 
\frac{\sqrt{acu}}{b} \cot^{-1} \left( \frac{2b}{\sqrt{acu}} \right).
\end{equation}
\end{Cor}

\medskip

The integral $\mathbf{3.252.2}$ 
\begin{equation}
\int_{-\infty}^{\infty} \frac{dx}{(ax^{2}+2bx+c)^{n}} = 
\frac{(2n-3)!! \pi \, a^{n-1}}{(2n-2)!! \, (ac-b^{2})^{n-1/2}}
\end{equation}
\noindent
reduces via $u = a(x+b/a)/\sqrt{ac-b^{2}}$ to Wallis' integral 
\begin{equation}
\ift \frac{du}{(u^{2}+1)^{n}} = \frac{(2n-3)!!}{(2n-2)!!} \, \frac{\pi}{2},
\end{equation}
\noindent
that appears as $\mathbf{3.249.1}$. The reader will find in \cite{irrbook} 
proofs of Wallis' integral. Observe that the evaluation of 
$\mathbf{3.252.2}$ is much simpler than the corresponding 
half-line example presented in Lemma \ref{ex2}. \\

The last example of this type is  $\mathbf{3.252.3}$:
\begin{equation}
\ift \frac{dx}{(ax^{2}+2bx+c)^{n+3/2}} = \frac{(-2)^{n}}{(2n+1)!!} 
\frac{\partial^{n}}{\partial c^{n}} \left( 
\frac{1}{\sqrt{c} \, ( \sqrt{ac} + b)} \right).  \label{three-halfs}
\nonumber
\end{equation}

\noindent
A simple trigonometric substitution gives the case $n=0$:
\begin{eqnarray}
\ift \frac{dx}{(ax^{2}+2bx+c)^{3/2}} & = & \frac{\sqrt{a}}{ac-b^{2}} 
\int_{b/\sqrt{-d}}^{\infty} \frac{du}{(u^{2}+1)^{3/2}} \nonumber \\
& = & \frac{1}{\sqrt{c} \, ( \sqrt{ac} + b) }. \nonumber
\end{eqnarray}
\noindent
The general case follows by differentiating with respect to $c$ and observing
that 
\begin{equation}
\left( \frac{\partial}{\partial c} \right)^{j} = 
(-1)^{j} \frac{(2j+1)!!}{2^{j}} (ax^{2}+bx+c)^{-3/2-j}. \nonumber
\end{equation}

We now provide a closed-form expression for (\ref{three-halfs}). 

\begin{Thm}
Let $a, \, b, \, c \in \mathbb{R}$ and $n \in \mathbb{N}$. Define
$u = (ac-b^{2})/4ac$ and assume $cu > 0$. Then
\begin{equation}
\ift \frac{dx}{(ax^{2}+2bx+c)^{n+3/2}}  = 
\frac{(cu)^{-n}}{\sqrt{c} \binom{2n}{n} (2n+1)} 
\left( \frac{1}{\sqrt{ac}+b} - \frac{b}{ac-b^{2}} 
\sum_{j=1}^{n} \binom{2j}{j}u^{j} \right).
\nonumber
\end{equation}
\end{Thm}
\begin{proof}
Change parameters sequentially as $a \mapsto a^{2}; \, c \mapsto c^{2}; \, 
b \mapsto abc$. Then, in the new format both sides 
satisfy the differential-difference equation 
\begin{equation}
-(2N(1-b^{2})c)f_{N+1} = \frac{df_{N}}{dc} - \frac{b}{ac^{2N}}, 
\end{equation}
\noindent
where $N = n + \tfrac{3}{2}$. 
\end{proof}

\section{An elementary combination of exponentials and rational functions} 
\label{sec-expoele} 
\setcounter{equation}{0}

The table \cite{gr} contains two integrals belonging to the family 
\begin{equation}
T_{j} := \ift e^{-px} (e^{-x}-1)^{n} \frac{dx}{x^{j}}. 
\end{equation}
\noindent
Indeed $\mathbf{3.411.19}$ gives $T_{1}$:
\begin{equation}
\ift e^{-px} (e^{-x}-1)^{n} \frac{dx}{x} = - \sum_{k=0}^{n} (-1)^{k} 
\binom{n}{k} \ln(p+n-k), 
\end{equation}
\noindent
and $\mathbf{3.411.20}$ gives $T_{2}$:
\begin{equation}
\ift e^{-px} (e^{-x}-1)^{n} \frac{dx}{x^{2}} = \sum_{k=0}^{n} (-1)^{k} 
\binom{n}{k} (p+n-k) \ln(p+n-k), 
\end{equation}

The next result presents an explicit evaluation of $T_{j}$.

\begin{Prop}
Let $p$ be a free parameter, and let $n, \, j \in \mathbb{N}$ with $n +p > 0$.
Then
\begin{equation}
\int_{0}^{\infty}  e^{-px} (e^{-x} -1)^{n} \frac{dx}{x^{j}} = 
\frac{(-1)^{j}}{(j-1)!} 
\sum_{k=0}^{n} (-1)^{k} (p+n-k)^{j-1} \ln(p+n-k). 
\end{equation}
\end{Prop}
\begin{proof}
Start with the observation that 
\begin{equation}
T_{j} = - \int T_{j-1}(p) dp + C. 
\end{equation}
\noindent
Therefore we need to describe the iterative integrals $f_{j}(p) = 
\int f_{j-1}(p) \, dp$,with $f_{0}(p) = \ln(p + \alpha)$. This can be found in 
page 82 of \cite{irrbook} as 
\begin{equation}
f_{j}(p) = \frac{1}{j!} (p + \alpha)^{j} \ln(p + \alpha) - 
\frac{H_{j}}{j!} (p+ \alpha)^{j} + C, 
\end{equation}
\noindent 
with $\alpha = p+n-k$ and $H_{j} = 1 + \tfrac{1}{2} + \cdots + \frac{1}{j}$ 
is the harmonic number. To build back the functions $T_{j}$ employ the 
fact that, for any polynomial $Q(n,k)$, 
\begin{equation}
\sum_{k=0}^{n} (-1)^{k} (-1)^{k} \binom{n}{k} Q(n,k) \equiv 0. 
\end{equation}
\noindent
Consequently,
\begin{equation}
T_{j} = C + \frac{(-1)^{j+1}}{j!} \sum_{k=0}^{n} (-1)^k (p+n-k)^{j} \ln(p+n-k).
\end{equation}
\noindent
The last step is to check that $C = 0$. This follows directly from $T_{j} \to 
0$ as $p \to \infty$. The assertion is now validated.
\end{proof}

\section{An elementary logarithmic integral} 
\label{sec-elelog2} 
\setcounter{equation}{0}

Entry $\mathbf{4.222.1}$ states that
\begin{equation}
\ift \ln \left( \frac{a^{2}+x^{2}}{b^{2}+x^{2}} \right) \, dx = (a-b) \pi.
\end{equation}
\noindent
In order to establish this, we consider the finite integral
\begin{equation}
I(m):= \int_{0}^{m}  
\ln \left( \frac{a^{2}+x^{2}}{b^{2}+x^{2}} \right) \, dx 
\end{equation}
\noindent
and then let $m \to \infty$. 

Integration by parts gives
\begin{eqnarray}
\int_{0}^{m} \ln(a^{2}+x^{2}) \, dx 
 & = & m \ln (m^{2}+a^{2}) - 2 \int_{0}^{m} \frac{x^{2} \, dx}{a^{2}+x^{2}}
\nonumber \\
     & = & m \ln (m^{2}+a^{2}) - 2m + 
        2a^{2} \int_{0}^{m} \frac{dx}{a^{2}+x^{2}} \nonumber \\
     & = & m \ln (m^{2}+a^{2}) - 2m +  2a \tan^{-1} \left( \frac{m}{a} 
\right). \nonumber 
\end{eqnarray}
\noindent
Therefore
\begin{equation}
I(m) = m \ln \left( \frac{m^{2}+a^{2}}{m^{2}+b^{2}} \right) 
 + 2a \tan^{-1} \left( \frac{m}{a} 
\right) - 2b \tan^{-1} \left( \frac{m}{b} \right).  \nonumber
\end{equation}
\noindent
The limit of the logarithmic part is zero and the arctangent part gives 
$(a-b) \pi$ as required. 

\medskip

The generalization
\begin{equation}
\int_{0}^{\infty} \ln \left( \frac{a^{s}+x^{s}}{b^{s} + x^{s}} \right) \, dx 
= (a-b) \frac{\pi}{\sin(\pi/s)}
\end{equation}
\noindent
can be established by elementary methods provided we assume the value 
\begin{equation}
\int_{0}^{\infty} \frac{dx}{1+x^{s}} = \frac{\pi}{s \, \sin(\pi/s)} 
\label{known}
\end{equation}
\noindent
as given. This integral is evaluated in terms of Euler's beta function 
in \cite{moll-gr6}. Indeed, integration by parts gives
\begin{equation}
\int_{0}^{y} \ln(a^{s}+x^{s}) \, dx = y \ln(a^{s}+y^{s}) - sy + sa^{s}
\int_{0}^{y} \frac{dx}{a^{s}+y^{s}}, 
\end{equation}
\noindent
and similarly for the $b$-parameter. Combining these evaluations gives
\begin{equation}
\int_{0}^{y} \ln \left(\frac{a^{s}+x^{s}}{b^{s}+x^{s}} \right) \, dx=
y \ln \left( \frac{a^{s}+y^{s}}{b^{s} + y^{s}} \right) + 
sa^{s} \int_{0}^{y} \frac{dx}{a^{s}+x^{s}} -
sb^{s} \int_{0}^{y} \frac{dx}{b^{s}+x^{s}}. 
\nonumber
\end{equation}
\noindent
Upon letting $y \to \infty$, we observe that the logarithmic term vanishes 
and a scaling reduces the remaining integrals to (\ref{known}). 

\bigskip

\noindent
{\bf Acknowledgments}. The second author acknowledges the partial support of 

\noindent
$\text{NSF-DMS } 0409968$.

\end{document}